\documentstyle[11pt]{article}
\title{The Product Formula for Gromov-Witten Invariants}
\author{K. Behrend}

\catcode`\@=11

\font\tenmsx=msam10
\font\sevenmsx=msam7
\font\fivemsx=msam5
\font\tenmsy=msbm10
\font\sevenmsy=msbm7
\font\fivemsy=msbm5
\newfam\msxfam
\newfam\msyfam
\textfont\msxfam=\tenmsx  \scriptfont\msxfam=\sevenmsx
  \scriptscriptfont\msxfam=\fivemsx
\textfont\msyfam=\tenmsy  \scriptfont\msyfam=\sevenmsy
  \scriptscriptfont\msyfam=\fivemsy

\def\hexnumber@#1{\ifcase#1 0\or1\or2\or3\or4\or5\or6\or7\or8\or9\or
	A\or B\or C\or D\or E\or F\fi }

\font\teneuf=eufm10
\font\seveneuf=eufm7
\font\fiveeuf=eufm5
\newfam\euffam
\textfont\euffam=\teneuf
\scriptfont\euffam=\seveneuf
\scriptscriptfont\euffam=\fiveeuf
\def\frak{\ifmmode\let\next\frak@\else
 \def\next{\Err@{Use \string\frak\space only in math mode}}\fi\next}
\def\goth{\relaxnext@\ifmmode\let\next\frak@\else
 \def\next{\Err@{Use \string\goth\space only in math mode}}\fi\next}
\def\frak@#1{{\frak@@{#1}}}
\def\frak@@#1{\fam\euffam#1}

\edef\msx@{\hexnumber@\msxfam}
\edef\msy@{\hexnumber@\msyfam}

\mathchardef\boxdot="2\msx@00
\mathchardef\boxplus="2\msx@01
\mathchardef\boxtimes="2\msx@02
\mathchardef\square="0\msx@03
\mathchardef\blacksquare="0\msx@04
\mathchardef\centerdot="2\msx@05
\mathchardef\lozenge="0\msx@06
\mathchardef\blacklozenge="0\msx@07
\mathchardef\circlearrowright="3\msx@08
\mathchardef\circlearrowleft="3\msx@09
\mathchardef\rightleftharpoons="3\msx@0A
\mathchardef\leftrightharpoons="3\msx@0B
\mathchardef\boxminus="2\msx@0C
\mathchardef\Vdash="3\msx@0D
\mathchardef\Vvdash="3\msx@0E
\mathchardef\vDash="3\msx@0F
\mathchardef\twoheadrightarrow="3\msx@10
\mathchardef\twoheadleftarrow="3\msx@11
\mathchardef\leftleftarrows="3\msx@12
\mathchardef\rightrightarrows="3\msx@13
\mathchardef\upuparrows="3\msx@14
\mathchardef\downdownarrows="3\msx@15
\mathchardef\upharpoonright="3\msx@16

\mathchardef\downharpoonright="3\msx@17
\mathchardef\upharpoonleft="3\msx@18
\mathchardef\downharpoonleft="3\msx@19
\mathchardef\rightarrowtail="3\msx@1A
\mathchardef\leftarrowtail="3\msx@1B
\mathchardef\leftrightarrows="3\msx@1C
\mathchardef\rightleftarrows="3\msx@1D
\mathchardef\Lsh="3\msx@1E
\mathchardef\Rsh="3\msx@1F
\mathchardef\rightsquigarrow="3\msx@20
\mathchardef\leftrightsquigarrow="3\msx@21
\mathchardef\looparrowleft="3\msx@22
\mathchardef\looparrowright="3\msx@23
\mathchardef\circeq="3\msx@24
\mathchardef\succsim="3\msx@25
\mathchardef\gtrsim="3\msx@26
\mathchardef\gtrapprox="3\msx@27
\mathchardef\multimap="3\msx@28
\mathchardef\therefore="3\msx@29
\mathchardef\because="3\msx@2A
\mathchardef\doteqdot="3\msx@2B

\mathchardef\triangleq="3\msx@2C
\mathchardef\precsim="3\msx@2D
\mathchardef\lesssim="3\msx@2E
\mathchardef\lessapprox="3\msx@2F
\mathchardef\eqslantless="3\msx@30
\mathchardef\eqslantgtr="3\msx@31
\mathchardef\curlyeqprec="3\msx@32
\mathchardef\curlyeqsucc="3\msx@33
\mathchardef\preccurlyeq="3\msx@34
\mathchardef\leqq="3\msx@35
\mathchardef\leqslant="3\msx@36
\mathchardef\lessgtr="3\msx@37
\mathchardef\backprime="0\msx@38
\mathchardef\risingdotseq="3\msx@3A
\mathchardef\fallingdotseq="3\msx@3B
\mathchardef\succcurlyeq="3\msx@3C
\mathchardef\geqq="3\msx@3D
\mathchardef\geqslant="3\msx@3E
\mathchardef\gtrless="3\msx@3F
\mathchardef\sqsubset="3\msx@40
\mathchardef\sqsupset="3\msx@41
\mathchardef\vartriangleright="3\msx@42
\mathchardef\vartriangleleft="3\msx@43
\mathchardef\trianglerighteq="3\msx@44
\mathchardef\trianglelefteq="3\msx@45
\mathchardef\bigstar="0\msx@46
\mathchardef\between="3\msx@47
\mathchardef\blacktriangledown="0\msx@48
\mathchardef\blacktriangleright="3\msx@49
\mathchardef\blacktriangleleft="3\msx@4A
\mathchardef\vartriangle="0\msx@4D
\mathchardef\blacktriangle="0\msx@4E
\mathchardef\triangledown="0\msx@4F
\mathchardef\eqcirc="3\msx@50
\mathchardef\lesseqgtr="3\msx@51
\mathchardef\gtreqless="3\msx@52
\mathchardef\lesseqqgtr="3\msx@53
\mathchardef\gtreqqless="3\msx@54
\mathchardef\Rrightarrow="3\msx@56
\mathchardef\Lleftarrow="3\msx@57
\mathchardef\veebar="2\msx@59
\mathchardef\barwedge="2\msx@5A
\mathchardef\doublebarwedge="2\msx@5B
\mathchardef\angle="0\msx@5C
\mathchardef\measuredangle="0\msx@5D
\mathchardef\sphericalangle="0\msx@5E
\mathchardef\varpropto="3\msx@5F
\mathchardef\smallsmile="3\msx@60
\mathchardef\smallfrown="3\msx@61
\mathchardef\Subset="3\msx@62
\mathchardef\Supset="3\msx@63
\mathchardef\Cup="2\msx@64

\mathchardef\Cap="2\msx@65

\mathchardef\curlywedge="2\msx@66
\mathchardef\curlyvee="2\msx@67
\mathchardef\leftthreetimes="2\msx@68
\mathchardef\rightthreetimes="2\msx@69
\mathchardef\subseteqq="3\msx@6A
\mathchardef\supseteqq="3\msx@6B
\mathchardef\bumpeq="3\msx@6C
\mathchardef\Bumpeq="3\msx@6D
\mathchardef\lll="3\msx@6E

\mathchardef\ggg="3\msx@6F

\mathchardef\circledS="0\msx@73
\mathchardef\pitchfork="3\msx@74
\mathchardef\dotplus="2\msx@75
\mathchardef\backsim="3\msx@76
\mathchardef\backsimeq="3\msx@77
\mathchardef\complement="0\msx@7B
\mathchardef\intercal="2\msx@7C
\mathchardef\circledcirc="2\msx@7D
\mathchardef\circledast="2\msx@7E
\mathchardef\circleddash="2\msx@7F
\def\ulcorner{\delimiter"4\msx@70\msx@70 }
\def\urcorner{\delimiter"5\msx@71\msx@71 }
\def\llcorner{\delimiter"4\msx@78\msx@78 }
\def\lrcorner{\delimiter"5\msx@79\msx@79 }
\def\yen{\mathhexbox\msx@55 }
\def\checkmark{\mathhexbox\msx@58 }
\def\circledR{\mathhexbox\msx@72 }
\def\maltese{\mathhexbox\msx@7A }
\mathchardef\lvertneqq="3\msy@00
\mathchardef\gvertneqq="3\msy@01
\mathchardef\nleq="3\msy@02
\mathchardef\ngeq="3\msy@03
\mathchardef\nless="3\msy@04
\mathchardef\ngtr="3\msy@05
\mathchardef\nprec="3\msy@06
\mathchardef\nsucc="3\msy@07
\mathchardef\lneqq="3\msy@08
\mathchardef\gneqq="3\msy@09
\mathchardef\nleqslant="3\msy@0A
\mathchardef\ngeqslant="3\msy@0B
\mathchardef\lneq="3\msy@0C
\mathchardef\gneq="3\msy@0D
\mathchardef\npreceq="3\msy@0E
\mathchardef\nsucceq="3\msy@0F
\mathchardef\precnsim="3\msy@10
\mathchardef\succnsim="3\msy@11
\mathchardef\lnsim="3\msy@12
\mathchardef\gnsim="3\msy@13
\mathchardef\nleqq="3\msy@14
\mathchardef\ngeqq="3\msy@15
\mathchardef\precneqq="3\msy@16
\mathchardef\succneqq="3\msy@17
\mathchardef\precnapprox="3\msy@18
\mathchardef\succnapprox="3\msy@19
\mathchardef\lnapprox="3\msy@1A
\mathchardef\gnapprox="3\msy@1B
\mathchardef\nsim="3\msy@1C
\mathchardef\ncong="3\msy@1D

\mathchardef\varsubsetneq="3\msy@20
\mathchardef\varsupsetneq="3\msy@21
\mathchardef\nsubseteqq="3\msy@22
\mathchardef\nsupseteqq="3\msy@23
\mathchardef\subsetneqq="3\msy@24
\mathchardef\supsetneqq="3\msy@25
\mathchardef\varsubsetneqq="3\msy@26
\mathchardef\varsupsetneqq="3\msy@27
\mathchardef\subsetneq="3\msy@28
\mathchardef\supsetneq="3\msy@29
\mathchardef\nsubseteq="3\msy@2A
\mathchardef\nsupseteq="3\msy@2B
\mathchardef\nparallel="3\msy@2C
\mathchardef\nmid="3\msy@2D
\mathchardef\nshortmid="3\msy@2E
\mathchardef\nshortparallel="3\msy@2F
\mathchardef\nvdash="3\msy@30
\mathchardef\nVdash="3\msy@31
\mathchardef\nvDash="3\msy@32
\mathchardef\nVDash="3\msy@33
\mathchardef\ntrianglerighteq="3\msy@34
\mathchardef\ntrianglelefteq="3\msy@35
\mathchardef\ntriangleleft="3\msy@36
\mathchardef\ntriangleright="3\msy@37
\mathchardef\nleftarrow="3\msy@38
\mathchardef\nrightarrow="3\msy@39
\mathchardef\nLeftarrow="3\msy@3A
\mathchardef\nRightarrow="3\msy@3B
\mathchardef\nLeftrightarrow="3\msy@3C
\mathchardef\nleftrightarrow="3\msy@3D
\mathchardef\divideontimes="2\msy@3E
\mathchardef\varnothing="0\msy@3F
\mathchardef\nexists="0\msy@40
\mathchardef\mho="0\msy@66
\mathchardef\eth="0\msy@67
\mathchardef\eqsim="3\msy@68
\mathchardef\beth="0\msy@69
\mathchardef\gimel="0\msy@6A
\mathchardef\daleth="0\msy@6B
\mathchardef\lessdot="3\msy@6C
\mathchardef\gtrdot="3\msy@6D
\mathchardef\ltimes="2\msy@6E
\mathchardef\rtimes="2\msy@6F
\mathchardef\shortmid="3\msy@70
\mathchardef\shortparallel="3\msy@71
\mathchardef\smallsetminus="2\msy@72
\mathchardef\thicksim="3\msy@73
\mathchardef\thickapprox="3\msy@74
\mathchardef\approxeq="3\msy@75
\mathchardef\succapprox="3\msy@76
\mathchardef\precapprox="3\msy@77
\mathchardef\curvearrowleft="3\msy@78
\mathchardef\curvearrowright="3\msy@79
\mathchardef\digamma="0\msy@7A
\mathchardef\varkappa="0\msy@7B
\mathchardef\hslash="0\msy@7D
\mathchardef\hbar="0\msy@7E
\mathchardef\backepsilon="3\msy@7F
\def\Bbb{\ifmmode\let\next\Bbb@\else
 \def\next{\errmessage{Use \string\Bbb\space only in math
 mode}}\fi\next} 
\def\Bbb@#1{{\Bbb@@{#1}}}
\def\Bbb@@#1{\fam\msyfam#1}

\catcode`\@=12

\newtheorem{prop}{Proposition}
\newtheorem{lem}[prop]{Lemma}

\newtheorem{them}[prop]{Theorem}

\newtheorem{defnp}[prop]{Definition}
\newtheorem{numconp}[prop]{Construction}
\newtheorem{numrmkp}[prop]{Remark}
\newtheorem{numexp}[prop]{Example}
\newtheorem{numrmksp}[prop]{Remarks}
\newtheorem{conp}[prop]{}

\newenvironment{numrmk}{\begin{numrmkp}\rm}{\end{numrmkp}}

\newtheorem{warningp}{Warning}
\newtheorem{notep}{Note}
\newtheorem{claimp}{Claim}
\newtheorem{examplep}{Example}
\newtheorem{examplesp}{Examples}
\newtheorem{rmkp}{Remark}
\newtheorem{rmksp}{Remarks}

\newenvironment{pf}{\begin{trivlist}\item[]{\sc Proof.}}%
            {\nolinebreak $\Box$ \end{trivlist}}

\newcommand{\noprint}[1]{}

\renewcommand{\tilde}{\widetilde}
\newcommand{\dual}[1]{#1^{\vee}}

\newcommand{\cart}{{\mbox{\tiny cart}}}

\newcommand{\stab}{\mbox{\tiny stab}}

\newcommand{\upst}{^{\ast}}

\newcommand{\upsh}{^{!}}
\newcommand{\lst}{_{\ast}}

\newcommand{\com}{^{\scriptscriptstyle\bullet}}

\newcommand{\DD}{{\frak D}}

\newcommand{\PP}{{\frak P}}

\newcommand{\GG}{{\frak G}}         

\newcommand{\MM}{{\frak M}}

\newcommand{\VV}{{\frak V}}

\newcommand{\qq}{{\Bbb Q}}
\newcommand{\pp}{{\Bbb P}}

\newcommand{\cC}{{\cal C}}

\newcommand{\del}{\partial}

\newcommand{\ob}{\mathop{\rm ob}}

\newcommand{\spec}{\mathop{\rm Spec}\nolimits}

\newcommand{\comp}{\mathbin{{\scriptstyle\circ}}}
\newcommand{\ol}{\overline}

\newcommand{\ldiag}[1]%
       {\makebox[0cm]{${\scriptstyle#1}\downarrow\phantom{\scriptstyle#1}$}}
\newcommand{\ldiagup}[1]%
       {\makebox[0cm]{${\scriptstyle#1}\uparrow\phantom{\scriptstyle#1}$}}
\newcommand{\rdiag}[1]%
       {\makebox[0cm]{$\phantom{\scriptstyle#1}\downarrow{\scriptstyle#1}$}}
\newcommand{\sediagr}[1]%
       {\makebox[0cm]{$\phantom{\scriptstyle#1}\searrow{\scriptstyle#1}$}}
\newcommand{\nediagr}[1]%
       {\makebox[0cm]{$\phantom{\scriptstyle#1}\nearrow{\scriptstyle#1}$}}
\newcommand{\rdiagup}[1]%
       {\makebox[0cm]{$\phantom{\scriptstyle#1}\uparrow{\scriptstyle#1}$}}
\newcommand{\swdiag}[1]%
       {\makebox[0cm]{$\phantom{\scriptstyle#1}\swarrow{\scriptstyle#1}$}}
\newcommand{\sediag}[1]%
       {\makebox[0cm]{${\scriptstyle#1}\searrow\phantom{\scriptstyle#1}$}}
\newcommand{\nediag}[1]%
       {\makebox[0cm]{${\scriptstyle#1}\nearrow\phantom{\scriptstyle#1}$}}

\newcommand{\longiso}{\stackrel{\textstyle\sim}{\longrightarrow}}

\newcommand{\doublearrowstack}[2]%
{{{{\scriptstyle#1}\atop{\textstyle\longrightarrow}}
\atop{{\textstyle\longrightarrow}\atop{\scriptstyle#2}}}}
\newcommand{\rightleftarrowstack}[2]%
{{{{\scriptstyle#1}\atop{\textstyle\longrightarrow}}
\atop{{\textstyle\longleftarrow}\atop{\scriptstyle#2}}}}
\newcommand{\leftrightarrowstack}[2]%
{{{{\scriptstyle#1}\atop{\textstyle\longleftarrow}}
\atop{{\textstyle\longrightarrow}\atop{\scriptstyle#2}}}}

\newcommand{\comdia}[9]{%
\begin{array}{ccc}
#1 & \stackrel{#2}{\longrightarrow} & #3 \\
\ldiag{#4} & #5 & \rdiag{#6} \\
#7 & \stackrel{#8}{\longrightarrow} & #9
\end{array}}

\newcommand{\overtoparrow}%
{\makebox[0cm]{\beginpicture
\setcoordinatesystem units <.8cm,.4cm> point at 0 0
\setplotarea x from -3 to 3, y from 0 to 1
\setquadratic
\plot -3 0 0 1 3 0 /
\put{\vector(3,-1){0}}[Bl] at 3 0
\endpicture}}

\newcommand{\underbottomarrow}%
{\makebox[0cm]{\beginpicture
\setcoordinatesystem units <.8cm,.4cm> point at 0 0
\setplotarea x from -3 to 3, y from 0 to 1
\setquadratic
\plot -3 1 0 0 3 1 /
\put{\vector(3,1){0}}[Bl] at 3 1
\endpicture}}

\renewcommand{\t}{_{\tau}}

\newcommand{\ses}[5]%
{0\longrightarrow#1\stackrel{#2}{ \longrightarrow}#3\stackrel{#4}{
\longrightarrow}#5\longrightarrow0}

\newcommand{\dt}[6]%
{#1\stackrel{#2}{longrightarrow}#3 \stackrel{#4}{\longrightarrow}#5
\stackrel{#6}{\longrightarrow} #1[1]}  
 
\newcommand{\cat}[1]%
{(\mbox{\rm #1})}

\begin{document} \sloppy
\date{October 10, 1997}
\maketitle

\begin{abstract}
We prove that the system of Gromov-Witten invariants of the product of
two varieties is equal to the tensor product of the systems of
Gromov-Witten invariants of the two factors.
\end{abstract}

\newcommand{\ms}{\ol{M}(W,\tau)}

\tableofcontents

\subsection{Introduction}

\newcommand{\gsv}{\mbox{$\tilde{\GG}_s(V)$}}
\newcommand{\gsw}{\mbox{$\tilde{\GG}_s(W)$}}
\newcommand{\gso}{\mbox{$\tilde{\GG}_s(0)$}}
\newcommand{\gsvw}{\mbox{$\tilde{\GG}_s(V\times W)$}}

Let $V$ and $W$ be smooth and projective varieties over the field
$k$. In this article we treat the question how to express the
Gromov-Witten invariants of $V\times W$ in terms of the Gromov-Witten
invariants of $V$ and $W$. 

On an intuitive level, the answer is quite obvious. For example,
assume $V=W=\pp^1$ and let us ask the question how many curves in
$\pp^1\times\pp^1$ of genus $g$ and bidegree $(d_1,d_2)$ pass through
$n=2(d_1+d_2)+g-1$ given points $P_1,\ldots,P_n$ of $\pp^1\times
\pp^1$ in general position. The answer is given by the Gromov-Witten
invariant $$I^{\pp^1\times\pp^1}_{g,n}(d_1,d_2)(\gamma^{\otimes n}),$$
where $\gamma\in H^4(\pp^1\times\pp^1,\qq)$ is the cohomology class
Poincar\'e dual to a point. 

We rephrase the question by asking how many triples
$(C,x_1,\ldots,x_n,f)$, where $C$ is a curve of genus $g$,
$x_1,\ldots,x_n$ are marked points on $C$ and $f:C\to\pp^1\times\pp^1$
is a morphism of bidegree $(d_1,d_2)$ exist (up to isomorphism) which
satisfy $f(x_i)=P_i$, for all $i=1,\ldots,n$. Now a morphism
$f:C\to\pp^1\times\pp^1$ of bidegree $(d_1,d_2)$ is given by two
morphisms $f_1:C\to\pp^1$ and $f_2:C\to\pp^1$ of degrees $d_1$ and
$d_2$, respectively. The requirement that $f(x_i)=P_i$ translates into
$f_1(x_i)=Q_i$ and $f_2(x_i)=R_i$, if we write the components of $P_i$
as $P_i=(Q_i,R_i)$. The family of all marked curves
$(C,x_1,\ldots,x_n)$ admitting such an $f_1$ is some cycle, say
$\Gamma_1$, in $\ol{M}_{g,n}$. Of course, the family of all curves
$(C,x_1,\ldots,x_2)$ admitting an $f_2$ as above is another cycle
$\Gamma_2$ in $\ol{M}_{g,n}$ and the family of all
$(C,x_1,\ldots,x_n)$ admitting an $f_1$ and an $f_2$ is the
intersection $\Gamma_1\cdot\Gamma_2$. So the Gromov-Witten number we
are interested in is
\[I^{\pp^1\times\pp^1}_{g,n}(d_1,d_2)(\gamma^{\otimes n}) =
\Gamma_1\cdot\Gamma_2.\]

In fact, the dual cohomology classes of $\Gamma_1$ and $\Gamma_2$ are
Gromov-Witten invariants themselves, namely
$I^{\pp^1}_{g,n}(d_1)(\tilde\gamma^{\otimes n})$ and
$I^{\pp^1}_{g,n}(d_2)(\tilde\gamma^{\otimes n})$, where
$\tilde\gamma\in H^2(\pp^1,\qq)$ is the cohomology class dual to a
point. Thus we have
\[I^{\pp^1\times\pp^1}_{g,n}(d_1,d_2)(\gamma^{\otimes n}) =
I^{\pp^1}_{g,n}(d_1)(\tilde\gamma^{\otimes n}) \cup
I^{\pp^1}_{g,n}(d_2)(\tilde\gamma^{\otimes n}) \] in
$H\upst(\ol{M}_{g,n}),\qq)$. This is the simplest instance of the
product formula, which we shall prove in this article. (Note that we
have identified, as usual, top degree cohomology classes on $\ol
M_{g,n}$ with their integrals over the fundamental cycle $[\ol
M_{g,n}]$.)

We get a more general statement by letting $V$ and $W$ be arbitrary
smooth projective varieties over $k$. We fix cohomology classes
$\gamma_1,\ldots,\gamma_n\in H\upst(V)$ and
$\epsilon_1,\ldots,\epsilon_n\in H\upst(W)$, which we assume to be
homogeneous, for simplicity. Then the product formula says that 
\begin{eqnarray}\label{bpf}
\lefteqn{
I^{V\times W}_{g,n}(\beta)
(\gamma_1\otimes\epsilon_1\otimes\ldots
\otimes\gamma_n\otimes\epsilon_n)}\nonumber\\
& = & (-1)^s
I^{V}_{g,n}(\beta_V)(\gamma_1\otimes\ldots \otimes\gamma_n) \cup 
I^{W}_{g,n}(\beta_W)(\epsilon_1\otimes\ldots \otimes\epsilon_n) 
\end{eqnarray}
in $H\upst(\ol{M}_{g,n},\qq)$. Here $\beta\in H_2(V\times W)^+$ and
$\beta_V={p_V}\lst\beta$, $\beta_W={p_W}\lst\beta$, where $p_V$ and
$p_W$ are the projections onto the factors of $V\times W$. The sign
is given by
\[s=\sum_{i>j}\deg\gamma_i\deg\epsilon_j.\]

This formula is already stated in \cite{KM} as a property expected of
Gromov-Witten invariants. In the case of $g=0$ and $V$ and $W$ (and
hence $V\times W$) convex, it is not difficult to prove, once the
properties of stacks of stable maps are established, as they are, for
example, in \cite{BM}. Essentially, the above intuitive argument can
then be translated into a rigorous proof. In the general case, the
enumerative meaning of Gromov-Witten invariants is much less clear,
since one has to use `virtual' fundamental classes to define
them. (This is done in \cite{BF} and \cite{gwi} or \cite{litian}.) So
the theorem follows from properties of virtual fundamental
classes. This is what we prove in the present paper.

Formula~(\ref{bpf}) has been used by various authors to understand the
quantum cohomology of a product. (See \cite{KMK}, \cite{KMZ} and
\cite{Kauf}.) By Formula~(\ref{bpf}), the codimension zero
Gromov-Witten invariants (ie.\ those that are numbers, like
$I^{\pp^1\times\pp^1}_{g,n}(d_1,d_2)(\gamma^{\otimes n})$, above) of a
product are determined by the Gromov-Witten invariants of higher
codimension of the factors and by the intersection theory of
$\ol{M}_{g,n}$.

To explain the treatment in this article, let us
reformulate~(\ref{bpf}) by saying that 
\[\begin{array}{ccc}
h(V\times W)^{\otimes n} & \stackrel{I^{V\times W}_{g,n}(\beta)}
{\longrightarrow} & h(\ol M_{g,n}) \\
\parallel & & \rdiagup{\Delta\upst} \\
h(V)^{\otimes n}\otimes h(W)^{\otimes n} &
\stackrel{I^{V}_{g,n}(\beta_V)\otimes
I^{W}_{g,n}(\beta_W)}{\longrightarrow} &  
h(\ol M_{g,n})\otimes h(\ol M_{g,n}),
\end{array}\]
where $\Delta:\ol{M}_{g,n}\to \ol{M}_{g,n}\times\ol{M}_{g,n}$ is the
diagonal, commutes. Here we have passed to the motivic Gromov-Witten
invariants. These are homomorphisms between DMC-motives. (These are
like Chow motives, except that they are made from smooth and proper
Deligne-Mumford stacks, instead of varieties. For details see
\cite{BM}, Section~8.)

To summarize all of their functorial properties, Gromov-Witten
invariants where defined in \cite{BM} as natural transformations
between the functors $h(V)^{\otimes S}$ and $h(\ol M)$, which are
functor from a certain graph category $\gsv_\cart$ to the category of
graded DMC-motives. To explain, let us start by reviewing some graph
theory. The category $\tilde{\GG}_s=\gso$ is the category of stable
modular graphs (graphs whose vertices are labeled with genuses; see
\cite{BM}, Definition~1.5) with so called extended isogenies as
morphisms. An {\em extended isogeny }is either a morphism gluing two
tails to an edge, or it is a proper isogeny (or a composition of the
two). An {\em isogeny }is a morphism which contracts various edges or
tails or both. (The name isogeny comes from the fact that such
morphisms do not affect the genus of the components of the graphs
involved.) For the definition of composition of extended isogenies,
see \cite{BM}, Page~36.

The category $\gsv_\cart$ is called the {\em cartesian extended
isogeny category }over $V$. The most important objects of $\gsv_\cart$
are pairs $(\tau,(\beta_i)_{i\in I})$, where $\tau$ is a stable
modular graph and $(\beta_i)_{i\in I}$ is a family of $H_2(V)^+$
markings on $\tau$. This means that each $\beta_i$ is a function
$\beta_i:V\t\to H_2(V)^+$, where $V\t$ is the set of vertices of
$\tau$. (The indexing set $I$ is finite.) The fundamental property of
$\gsv_\cart$ is that it is {\em fibered } over $\tilde{\GG}_s$. This
means that there is a functor $\gsv_\cart\to\tilde{\GG}_s$ (projection
onto the first component) and that given an object
$(\tau,(\beta_i)_{i\in I})$ of $\gsv_\cart$ and a morphism
$\phi:\sigma\to\tau$ there exists, up to isomorphism, a unique object
$(\sigma,(\gamma_j)_{j\in J})$ of $\gsv_\cart$ together with a
morphism $\Phi:(\sigma,(\gamma_j)_{j\in J}):\to(\tau,(\beta_i)_{i\in
I})$ covering $\phi:\sigma\to\tau$. When constructing $\Phi$, the
basic non-obvious case is that where $\phi$ contracts a non-looping
edge of $\sigma$ and $I$ has only one element. Then we have the graph
$\tau$ with an $H_2(V)^+$-marking $\beta$ and there are two vertices
$v_1$, $v_2$ of $\sigma$ corresponding to one vertex $w$ of
$\tau$. Then $(\sigma,(\gamma_j)_{j\in J})$ is defined such that $J$
counts the ways to write $\beta(w)=\beta_1+\beta_2$ in $H_2(V)^+$ and
$\gamma_j$ assigns $\beta_1$ to $v_1$ and $\beta_2$ to $v_2$, and
otherwise does not differ from $\beta$. 

Things get more complicated, if one also considers the less important
objects of $\gsv_\cart$. These are of the form $(\tau,(\tau_i)_{i\in
I})$, where, as above, $\tau$ is a stable modular graph, but now each
$\tau_i$ is a stable $H_2(V)^+$-marked graph (as opposed to an
$H_2(V)^+$-marked stable graph), together with a {\em stabilizing
morphism }$\tau_i\to\tau$. For the complete picture, see \cite{BM},
Definition~5.9. 

On objects, the morphisms $h(V)^{\otimes S}$ and $h(\ol M)$ from
$\gsv_\cart$ to $(\mbox{graded DMC-motives})$ are defined as follows: 
For an object $(\tau,(\beta_i)_{i\in I})$ of $\gsv_\cart$ we have
\[h(V)^{\otimes S}(\tau,(\beta_i)_{i\in I})=h(V)^{\otimes S\t},\]
where $S\t$ is the set of tails of $\tau$ and 
\[h(\ol M)(\tau,(\beta_i)_{i\in I})=h(\ol M(\tau)),\]
where $$\ol M(\tau)=\prod_{v\in V\t}\ol{M}_{g(v),F\t(v)}$$ and
$F\t(v)$ is the set of flags meeting the vertex $v$ of $\tau$. So both
of these functors only depend on the first component $\tau$ of
$(\tau,(\beta_i)_{i\in I})$. For the definition of these functors on
morphisms, see \cite{BM}, Section~9. Note that $h(V)^{\otimes S}$
actually comes with a twist (ie.\ a degree shift) $\chi\dim V$. This
we shall ignore here, to shorten notation and since nothing
interesting happens to it, anyway. 

The {\em Gromov-Witten transformation } of $V$ is now defined as a
natural transformation
\[I^V:h(V)^{\otimes S}\longrightarrow h(\ol M)\]
of functors from $\gsv_\cart$ to $(\mbox{graded DMC-motives})$. In
this paper (Theorem~\ref{st}), we shall prove that 
\[I^{V\times W}=I^V\cup I^W,\]
where $I^V\cup I^W$ is defined as $\Delta\upst(I^V\otimes I^W)$. 

Since Gromov-Witten invariants are defined in terms of virtual
fundamental classes on moduli stacks of stable maps, this theorem
follows from a certain compatibility between virtual fundamental
classes. This is our main result (Theorem~\ref{pt}) and takes up most
of this paper.

\subsection{Virtual Fundamental Classes}

Fix a ground field $k$. For a smooth projective $k$-variety $V$ let
$\gsv$ be the category of extended isogenies of stable
$H_2(V)^+$-graphs bounded by the characteristic of $k$ (see \cite{BM},
Definition~5.6 and Example~II following Definition~5.11).  Let
$J(V,\tau)\in A_{\dim(V,\tau)}(\ol{M}(V,\tau))$, for $\tau\in\ob\gsv$,
be the `virtual fundamental class', or orientation (\cite{BM},
Definition~7.1) of $\ol{M}$ over $\gsv$ constructed in \cite{gwi},
Theorem~6, using the techniques from \cite{BF}.

Now let us consider two smooth projective $k$-varieties $V$ and $W$;
denote the two projections by $p_V:V\times W\to V$ and $p_W:V\times
W\to W$. If $\tau$ is a stable $H_2(V\times W)^+$-graph, we denote by
${p_V}\lst(\tau)$ and ${p_W}\lst(\tau)$ the stabilizations of $\tau$
with respect to ${p_V}\lst:H_2(V\times W)^+\to H_2(V)^+$ and
${p_W}\lst:H_2(V\times W)^+\to H_2(W)^+$ (see \cite{BM}, Remark~1.15),
by $\tau^s$ the absolute stabilization of $\tau$.

Applying the functor $\ol{M}$ to the commutative diagram
\[\comdia{(V\times
W,\tau)}{}{(W,{p_W}\lst(\tau))}{}{}{}{(V,{p_V}\lst(\tau))}{}{(\spec
k,\tau^s)} \]
in $\VV\GG_s$ (see \cite{BM}, Remark~3.1 and the remark following
Theorem~3.14) we get a commutative diagram of proper Deligne-Mumford
stacks 
\[\comdia{\ol{M}(V\times W,\tau)}{}{\ol{M}(W,{p_W}\lst(\tau))} {}{}{}
{\ol{M}(V,{p_V}\lst(\tau))}{}{\ol{M}(\tau^s).}\]
In general, this diagram is not cartesian; let $P$ be the cartesian
product
\[\comdia{P}{}{\ol{M}(W,{p_W}\lst(\tau))} {}{}{}
{\ol{M}(V,{p_V}\lst(\tau))}{}{\ol{M}(\tau^s).}\]
Rewrite these diagrams as follows:
\[\begin{array}{ccccc}
\ol{M}(V\times W,\tau) & \stackrel{h}{\longrightarrow} & P & 
\longrightarrow & \ol{M}(V,{p_V}\lst\tau)\times\ol{M}(W,{p_W}\lst\tau)
\\ & \searrow & \ldiag{} & & \rdiag{} \\
& & \ol{M}(\tau^s) & \stackrel{\Delta}{\longrightarrow} & 
\ol{M}(\tau^s)\times\ol{M}(\tau^s).\end{array}\]
To shorten notation, write $J(V\times W)=J(V\times W,\tau)$,
$J(V)=J(V,{p_V}\lst\tau)$ and $J(W)=J(W,{p_W}\lst\tau)$. 

\begin{them} \label{pt}
We have 
\[\Delta\upsh(J(V)\times J(W))=h\lst(J(V\times W)).\]
\end{them}

For a stable $A$-graph $\tau$ ($A=H_2(V\times W)^+$, $H_2(V)^+$ etc.)
we denote by $\MM(\tau)$ the algebraic $k$-stack of $\tau$-marked
prestable curves, forgetting the $A$-structure, and thinking of $\tau$
simply as a (possibly not stable) modular graph. We consider the
diagram
\begin{equation} \label{bigd} \begin{array}{cccccc}
& \ol{M}(V\times W,\tau) & \stackrel{h}{\longrightarrow} & P &
\longrightarrow & \ol{M}(V,{p_V}\lst\tau)\times\ol{M}(W,{p_W}\lst\tau)
\\ \phantom{M}/ & \rdiag{c} & & \rdiag{} & & \rdiag{a} \\ b\mid &
\DD(\tau) & \stackrel{l}{\longrightarrow} & \PP &
\stackrel{\phi}{\longrightarrow} &
\MM({p_V}\lst\tau)\times\MM({p_W}\lst\tau) \\ \phantom{nM}\searrow&
\rdiag{e} & \searrow & \rdiag{} & & \rdiag{s\times s} \\ & \MM(\tau) &
& \ol{M}(\tau^s) & \stackrel{\Delta}{\longrightarrow} &
\ol{M}(\tau^s)\times\ol{M}(\tau^s). \end{array}\end{equation} Here
$s\times s$ is given by stabilizations and $\PP$ is defined as the
fibered product of $\Delta$ and $s\times s$. The morphisms $a$ and $b$
are given by forgetting maps, retaining only marked curves.

The algebraic stack $\DD(\tau)$ is defined as follows. 
For a $k$-scheme $T$ the groupoid $\DD(\tau)(T)$ has as objects
diagrams
\begin{equation}\label{ddd}
\begin{array}{ccc}
(C,x) & \longrightarrow & (C'',x'') \\
\ldiag{} & & \\
(C',x') & & \end{array}\end{equation}
where $(C,x)$ is a $\tau$-marked prestable curve over $T$, $(C',x')$ a
${p_V}\lst(\tau)$-marked prestable curve over $T$ and $(C'',x'')$ a
${p_W}\lst(\tau)$-marked prestable curve over $T$. The arrow
$(C,x)\to(C',x')$ is a morphism of marked prestable curves covering
the morphism $\tau\to{p_V}\lst(\tau)$ of modular graphs. Similarly, 
$(C,x)\to(C'',x'')$ is a morphism of marked prestable curves covering
$\tau\to{p_W}\lst(\tau)$. This concept has not been defined in
\cite{BM}; the definition (in this special case) is as follows. Let us
explain it for the case of $W$ instead of $V$, since this will lead to
less confusion of notation with the set of vertices of a graph. The
morphism $\tau\to{p_W}\lst(\tau)$ is given by a combinatorial morphism
of $0$-marked graphs $a:{p_W}\lst(\tau)\to\tau$ (see \cite{BM},
Definition~1.7). So there are maps $a:V_{{p_W}\lst(\tau)}\to V_{\tau}$
and $a:F_{{p_W}\lst(\tau)}\to F_{\tau}$. The morphism
$(C,x)\to(C'',x'')$ is given by a family
$p=(p_v)_{v\in V_{{p_W}\lst(\tau)}}$ of morphisms of prestable curves
(\cite{BM}, Definition~2.1) $p_v:C_{a(v)}\to C''_v$ such that for
every $i\in F_{{p_W}\lst(\tau)}$ we have
$p_{\del(i)}(x_{a(i)})=x''_i$. 

\newcommand{\pv}{{p_V}}
\newcommand{\pw}{{p_W}}

There are morphisms of stacks $e:\DD(\tau)\to\MM(\tau)$,
$\DD(\tau)\to\MM(\pv\lst\tau)$ and $\DD(\tau)\to\MM(\pw\lst\tau)$,
given, respectively, by mapping Diagram~(\ref{ddd}) to $(C,x)$,
$(C',x')$ and $(C'',x'')$. Let us denote the product of the latter two
by  
\[\tilde{\Delta}:\DD(\tau)\longrightarrow
\MM(\pv\lst\tau)\times\MM(\pw\lst\tau). \]

\begin{lem} \label{lem2}
In Diagram~(\ref{ddd}) both morphisms induce isomorphisms on
stabilizations. 
\end{lem}
\begin{pf}
This follows from the fact that any morphism of stable marked curves
(with identical dual graphs) is an isomorphism. This fact is proved
in \cite{BM}, at the very end of the proof of Theorem~3.6, which
immediately precedes Definition~3.13.
\end{pf}

By this lemma there is a commutative diagram
\[\begin{array}{ccc}
\DD(\tau) & \stackrel{\tilde{\Delta}}{\longrightarrow} &
\MM(\pv\lst\tau)\times \MM(\pw\lst\tau) \\
\ldiag{e} & & \\
\MM(\tau) & & \rdiag{s\times s} \\
\ldiag{s} & & \\
\ol{M}(\tau^s) & \stackrel{\Delta}{\longrightarrow} &
\ol{M}(\tau^s)\times\ol{M}(\tau^s), \end{array}\]
which gives rise to the morphism $l:\DD(\tau)\to\PP$ of
Diagram~(\ref{bigd}). 

\begin{prop} \label{podd}
The morphisms $\Delta$ and $\tilde{\Delta}$ are proper regular local
immersions. Their natural orientations satisfy
$$l\lst[\tilde{\Delta}]=(s\times s)\upst[\Delta].$$
\end{prop}
\begin{pf}
Let $S_1$ and $S_2$ be finite sets, set $S=S_1\amalg S_2$. Let the
modular graph $\pv\lst(\tau)'$ be obtained from $\pv\lst(\tau)$ by
adding (in any fashion) $S_1$ to the set of tails of
$\pv\lst(\tau)$. Similarly, let $\pw\lst(\tau)'$ be obtained form
$\pw\lst(\tau)$ by adding $S_2$ to the set of tails, arbitrarily. Now
let $\tau'$ be obtained from $\tau$ by adding the set $S$ to the tails
of $\tau$ in the unique way such that $\tau\to\pv\lst(\tau)$ induces a
morphism $\tau'\to\pv\lst(\tau)'$, which gives the inclusion
$S_1\subset S$ on tails, and $\tau\to\pw\lst(\tau)$ induces a morphism
$\tau'\to\pw\lst(\tau)'$, which gives the inclusion $S_2\subset S$ in
tails. 

With these choices we have a cartesian diagram of $k$-stacks
\begin{equation}\label{poddi}
\comdia{\ol{M}(\tau')}{\delta} {\ol{M}(\pv\lst(\tau)') \times
\ol{M}(\pw\lst(\tau)')} {}{}{\chi} {\DD(\tau)} {\tilde{\Delta}}
{\MM(\pv\lst(\tau)) \times \MM(\pw\lst(\tau)).}
\end{equation}
The proof that this is the case is similar to the proof of
Proposition~\ref{prop5}, below. The morphism $\chi$ in
Diagram~(\ref{poddi}) is a local presentation of $\MM(\pv\lst(\tau))
\times \MM(\pw\lst(\tau))$, (see \cite{gwi}, remarks following
Lemma~1). Moreover, by choosing $S$ and the primed graphs correctly,
any point of $\MM(\pv\lst(\tau)) \times \MM(\pw\lst(\tau))$ can be
assumed to be in the image of $\chi$. So to prove that
$\tilde{\Delta}$ is a proper regular local immersion, it suffices to
prove that $\delta$ is a proper regular local immersion. Properness is
clear; the stacks $\ol{M}(\tau')$, $\ol{M}(\pv\lst(\tau)')$ and
$\ol{M}(\pw\lst(\tau)')$ are proper. The regular local immersion
property follows from injectivity on tangent spaces which can be
proved by a deformation theory argument. 

The proof for $\Delta$ is comparatively trivial. 

To prove the fact about the orientations, first note that $s\times s$
is flat (see \cite{gwi}, Proposition~3) and so $\phi$ is a regular local
immersion and $(s\times s)\upst[\Delta]=[\phi]$. To prove that
$l\lst[\tilde{\Delta}]=[\phi]$, it suffices to identify dense open
substacks $\DD(\tau)'\subset\DD(\tau)$ and $\PP'\subset\PP$ such that
$l$ induces an isomorphism $\DD(\tau)'\to\PP'$. We define $\DD(\tau)'$
to be the open substack of $\DD(\tau)$ over which $C_v\to C_v'$ is an
isomorphism for all $v\in V_{\pv\lst(\tau)}$ and $C_v\to C_v''$ is an
isomorphism for all $v\in V_{\pw\lst(\tau)}$. We define $\PP'$ to be
the pullback via $\Delta$ of $\MM(\pv\lst(\tau))' \times
\MM(\pw\lst(\tau))'$, where $\MM(\pv\lst(\tau))'$ is the open substack
over which $(C_v,(x_i)_{i\in F(v)})$ is stable, for all $v\in
V_{\tau^s}$, similarly for $\MM(\pw\lst(\tau))'$. Note the slight
abuse of notation; we have denoted vertices of different graphs by the
same letter. 
\end{pf}

\begin{lem}
The morphism $e:\DD(\tau)\to\MM(\tau)$ is \'etale.
\end{lem}
\begin{pf}
Similar to the proof of \cite{gwi}, Lemma~7.
\end{pf}

To complete Diagram~(\ref{bigd}), define a morphism $\ol{M}(V\times
W,\tau)\to\DD(\tau)$ by mapping a stable $(V\times W,\tau)$-map
$(C,x,f)$ first to the diagram
\[\begin{array}{ccc}
(C,x,f) & \longrightarrow & (C,x,\pw\comp f)^{\stab} \\
\ldiag{} & & \\
(C,x,\pv\comp f)^{\stab} & & \end{array}\]
and then passing to the underlying prestable curves. 

\begin{prop}\label{prop5}
The diagram
\[\comdia{\ol{M}(V\times W,\tau)} {}
{\ol{M}(V,\pv\lst\tau)\times\ol{M}(W,\pw\lst\tau)} {c}{}{a}{\DD(\tau)}
{\tilde{\Delta}} {\MM(\pv\lst\tau)\times\MM(\pw\lst\tau)}\]
is cartesian. 
\end{prop}
\begin{pf}
We have to construct a morphism from the fibered product of
$\tilde{\Delta}$ and $a$ to $\ol{M}(V\times W,\tau)$. So let there be
given a Diagram~(\ref{ddd}), representing an object of $\DD(\tau)(T)$,
for a $k$-scheme $T$. Moreover, let there be given families of maps
$(f')_{v\in V_{\pv\lst\tau}}$, $f_v':C_v'\to V$ and $(f'')_{v\in
V_{\pw\lst\tau}}$, $f_v'':C_v''\to W$, making $(C',x',f')$ and
$(C'',x'',f'')$ stable maps.  We need to construct a stable map from
$(C,x)$ to $V\times W$.  So let $v\in V_{\tau}$ be a vertex of
$\tau$. 

Let us construct a map $h_v:C_v\to W$. In case $v$ is in the
image of $V_{\pw\lst\tau}\to V_{\tau}$, and $w\mapsto v$ under this
map, we take $h_v$ to be the composition
$$C_v\stackrel{\pw}{\longrightarrow}
C_w''\stackrel{f_w''}{\longrightarrow} W.$$ In case $v$ is not in the
image of $V_{\pw\lst\tau}\to V_{\tau}$, then $v$ partakes in a long
edge or a long tail associated to and edge $\{i,\ol{i}\}$ or a tail
$i$ of $\pw\lst\tau$ (see the discussion of {\em stabilizing
morphisms}, Definition~5.7, in \cite{BM} for this terminology). Then
we define $f_v:C_v\to W$ to be the composition
$$C_v\longrightarrow T\stackrel{x_i''}{\longrightarrow}
C_{\del(i)}''\stackrel{f_{\del(i)}''}{\longrightarrow} W.$$ 
In the same manner, construct a map $g_v:C_v\to V$. Finally, let
$f_v:C_v\to V\times W$ be the product $g_v\times h_v$. Then the family
$(f_v)_{v\in V_{\tau}}$ makes $(C,x,f)$ a stable map over $T$ to
$V\times W$. One checks that $(C,x,\pv\comp f)^{\stab}=(C',x',f')$ and
$(C,x,\pw\comp f)^{\stab}=(C'',x'',f'')$, using the universal mapping
property of stabilization and the fact already alluded to in the proof
of Lemma~\ref{lem2}.
\end{pf}

Let $E\com(V)=E\com(V,\pv\lst\tau)$ and
$E\com(W)=E\com(W,\pw\lst\tau)$ denote the relative obstruction
theories for $\ol{M}(V,\pv\lst\tau)\to\MM(\pv\lst\tau)$ and
$\ol{M}(W,\pw\lst\tau)\to\MM(\pw\lst\tau)$, respectively, which were
defined in \cite{gwi}. As in \cite{BF} Proposition~7.4 there is an
induced obstruction theory $E\com(V)\boxplus E\com(W)$ for the
morphism $a$. Pulling back via $\tilde{\Delta}$ (as in \cite{BF}
Proposition~7.1) we get an induced obstruction theory
$\tilde{\Delta}\upst(E\com(V)\boxplus E\com(W))$ for the morphism
$c$. 

On the other hand, we have the relative obstruction theory
$E\com(V\times W)=E\com(V\times W,\tau)$ for the morphism $b$. Since
$e:\DD(\tau)\to\MM(\tau)$ is \'etale, we may think of $E\com(V\times
W)$ as a relative obstruction theory for $c$. 

\begin{prop}
The two relative obstruction theories
$\tilde{\Delta}\upst(E\com(V)\boxplus E\com(W))$ and $E\com(V\times
W)$ for the morphism $c$ are naturally isomorphic. 
\end{prop}
\begin{pf}
Let $\cC(V\times W,\tau)\to\ol{M}(V\times W,\tau)$,
$\cC(V,\pv\lst\tau)\to\ol{M}(V,\pv\lst\tau)$ and
$\cC(W,\pw\lst\tau)\to\ol{M}(W,\pw\lst\tau)$ be the universal
curves. Recall from \cite{gwi} that they are constructed by gluing the
curves associated to the vertices of a graph according the the edges
of that graph. Let us denote the pullbacks of the latter two universal
curves to $\ol{M}(V\times W,\tau)$ by $\cC(V)$ and $\cC(W)$,
respectively. We have maps $f_{V\times W}$, $f_V$ and $f_W$,
constructed from the universal stable maps, which fit into the
following commutative diagram:
\[\begin{array}{ccccc}
V & \stackrel{\pv}{\longleftarrow} & V\times W &
\stackrel{\pw}{\longrightarrow} & W \\
\ldiagup{f_V} & & \rdiagup{f_{V\times W}} & & \rdiagup{f_W} \\
\cC(V) & \stackrel{q_{V}}{\longleftarrow} & \cC(V\times W,\tau) &
\stackrel{q_W}{\longrightarrow} & \cC(W) \\
& \sediag{\pi_V} & \rdiag{\pi_{V\times W}} & \swdiag{\pi_W} & \\
&& \ol{M}(V\times W,\tau) & & \end{array}\]
By base change it is clear that 
$$\tilde{\Delta}\upst(\dual{E\com(V)}\boxplus\dual{E\com(W)})
= R{\pi_V}\lst f_V\upst T_V\oplus R{\pi_W}\lst f_W\upst T_W.$$
For any vector bundle $F$ on $\cC(W)$ the canonical homomorphism $F\to
R{q_W}\lst q_W\upst F$ is an isomorphism. Of course, the same property
is enjoyed by $q_V$. Hence we have a canonical isomorphism 
\begin{eqnarray*}
\lefteqn{R{\pi_V}\lst f_V\upst T_V\oplus R{\pi_W}\lst f_W\upst T_W
\longiso }\\
&&R{\pi_V}\lst R{q_V}\lst q_V\upst f_V\upst T_V\oplus R{\pi_W}\lst
R{q_W}\lst q_W\upst f_W\upst T_W \\
& = & R{\pi_{V\times W}}\lst f_{V\times W}\upst(T_V\boxplus T_W) \\
& = & \dual{E\com(V\times W)}.
\end{eqnarray*}
To conclude, we have a canonical isomorphism 
$$\tilde{\Delta}\upst(\dual{E\com(V)}\boxplus\dual{E\com(W)}) \longiso
\dual{E\com(V\times W)}$$
and by dualizing
$$E\com(V\times W) \longiso
\tilde{\Delta}\upst({E\com(V)}\boxplus{E\com(W)}).$$
\end{pf}

By this proposition, we have
\begin{eqnarray*}
J(V\times W) & = & [\ol{M}(V\times W,\tau), E\com(V\times W)] \\
             & = & [\ol{M}(V\times W,\tau),
             \tilde{\Delta}\upst(E\com(V)\boxplus E\com(W))] \\ 
             & = & \tilde{\Delta}\upsh [ \ol{M}(V,\pv\lst\tau) \times
             \ol{M}(W,\pw\lst\tau) , E\com(V)\boxplus E\com(W)] \\  
             &   & \quad\quad\mbox{(by \cite{BF} Proposition~7.2)} \\
             & = & \tilde{\Delta}\upsh (
             [\ol{M}(V,\pv\lst\tau),E\com(V)]  \times
             [\ol{M}(W,\pw\lst\tau),E\com(W)]) \\ 
             &   & \quad\quad\mbox{(by \cite{BF} Proposition~7.4)} \\       
             & = & \tilde{\Delta}\upsh (J(V)\times J(W)).
\end{eqnarray*}
So we may now calculate as follows:
\begin{eqnarray*}
\tilde{\Delta}\upsh (J(V)\times J(W)) & = & a\upst(s\times s)\upst
[\Delta] (J(V)\times J(W)) \\
   & = & a\upst l\lst[\tilde{\Delta}] (J(V)\times J(W)) \\
   &   & \quad\quad\mbox{(by Proposition~\ref{podd})} \\
   & = & h\lst \tilde{\Delta}\upsh(J(V)\times J(W)) \\
   & = & h\lst J(V\times W),
\end{eqnarray*}
which is the product property. This finishes the proof of
Theorem~\ref{pt}. 

\subsection{Gromov-Witten Transformations}

Theorem~\ref{pt} easily implies that the system of Gromov-Witten
invariants for $V\times W$ is equal to the tensor product (see
\cite{KM}, 2.5) of the systems of Gromov-Witten invariants for $V$ and
$W$, respectively. To get the full `operadic' picture, we need a few
graph theoretic preparations. 

\begin{prop} \label{psi}
There is a natural functor of categories fibered over $\gso$
\[\Psi:\gsvw_{\cart} \longrightarrow
\gsv_{\cart}\times_{\gso}\gsw_{\cart}.\] 
This functor is cartesian. 
\end{prop}
\begin{pf}
Let $p:V\to W$ be a morphism of smooth projective varieties over
$k$. We shall construct a natural functor of fibered categories over
$\gso$
$$\Psi_p:\gsv_{\cart}\longrightarrow\gsw_{\cart}.$$
This functor will be cartesian. 

For an object $(\tau,(\ol{a}_i,\tau_i)_{i\in I})$ of $\gsv_{\cart}$
let the image under $\Psi_p$ be $(\tau,(\ol{b}_i,p\lst(\tau_i))_{i\in
I})$. Here $p\lst(\tau_i)$ is the stabilization of $\tau_i$ covering
$p\lst:H_2(V)^+\to H_2(W)^+$. It comes with a natural morphism
$\ol{b}_i:p\lst(\tau_i)\to \tau$. To make $\ol{b}_i$ a stabilizing
morphism, we have to endow it with an orbit map (see \cite{BM},
Definition~5.7). Let $\ol{a}_i^m:E\t\cup S\t\to E_{\tau_i}\cup
S_{\tau_i}$ be the orbit map of $\ol{a}_i:\tau_i\to \tau$. Let
$\{f,\ol{f}\}$ be an edge of $\tau$. Then there exists a unique factor
$\{f',\ol{f}'\}$ of the long edge of $p\lst(\tau_i)$ associated to
$\{f,\ol{f}\}$, such that $\ol{a}_i^m(\{f,\ol{f}\})$ is a factor of
the long edge of $\tau_i$ associated to $\{f',\ol{f}'\}$. We set
$\ol{b}_i^m(\{f,\ol{f}\})=\{f',\ol{f}'\}$. This defines the orbit map
$\ol{b}_i^m$ of $\ol{b}_i$ on edges. On tails it is defined entirely
analogously. 

This achieves the definition of $\Psi_p$ on objects. We leave it to
the reader to explicate the action of $\Psi_p$ on morphisms; it boils
down to checking that the {\em pullback }of \cite{BM}, Definition~5.8
is compatible with applying $p\lst$. 

Now we define $\Psi$ by 
\begin{eqnarray*}
\Psi:\gsvw_{\cart}& \longrightarrow &
\gsv_{\cart}\times_{\gso}\gsw_{\cart} \\
(\tau,(\tau_i)) & \longmapsto &
(\Psi_{\pv}(\tau,(\tau_i)),\Psi_{\pw}(\tau,(\tau_i))). 
\end{eqnarray*}
\end{pf}

Let us denote, for any $V$, the Gromov-Witten transformation for $V$
(see \cite{BM}, Theorem~9.2), by
\[I^V:h(V)^{\otimes S}(\chi\dim V)\longrightarrow h(\ol{M}).\]
Recall that $I^V$ is a natural transformation between functors
\[\gsv_{\cart}\longrightarrow(\mbox{graded DMC-motives}),\]
where the two functors $h(V)^{\otimes S}(\chi\dim V)$ and $h(\ol{M})$
are induced from functors (with the same names)
$\gso\to(\mbox{DMC-motives})$, constructed in \cite{BM}, Section~9. 

\begin{numrmk} \label{lem8}
By our various definitions we have
\[h(V)^{\otimes S}(\chi\dim V)\otimes h(W)^{\otimes S}(\chi\dim W) =
h(V\times W)^{\otimes S}(\chi\dim V\times W)\]
as functors $\gso\to(\mbox{\rm DMC-motives})$. 
\end{numrmk}

The transformations $I^V$ and $I^W$ induce a transformation
\[I^V\otimes I^W: h(V)^{\otimes S}(\chi\dim V) \otimes h(W)^{\otimes
S} (\chi\dim W) \longrightarrow h(\ol{M})\otimes h(\ol{M})\]
between functors 
\[\gsv_{\cart}\times_{\gso}\gsw_{\cart}\longrightarrow (\mbox{graded
DMC-motives}).\]
It is defined as follows. Let $((\tau,(\tau_i)_{i\in
I}),(\tau,(\sigma_j)_{j\in J}))$ be an object of
$\gsv_{\cart}\times_{\gso}\gsw_{\cart}$. The value of $I^V\otimes I^W$
on this object is the morphism 
\[I^V(\tau,(\tau_i))\otimes I^W(\tau,(\sigma_j)): \]\[
h(V)^{\otimes 
S\t}(\chi(\tau)\dim V)\otimes h(W)^{\otimes S\t}(\chi(\tau)\dim W)
\longrightarrow h(\ol{M}(\tau))\otimes h(\ol{M}(\tau)).\]

Composing with $\Delta\upst:h(\ol{M})\otimes h(\ol{M})\to h(\ol{M})$
we get the transformation 
$$\Delta\upst(I^V\otimes I^W):h(V)^{\otimes S}(\chi\dim V) \otimes
h(W)^{\otimes S} (\chi\dim W) \longrightarrow h(\ol{M}),$$ 
which we shall also denote by $I^V\cup I^W=\Delta\upst(I^V\otimes
I^W)$. 

Pulling back via the functor $\Psi$ of Proposition~\ref{psi} and using
Remark~\ref{lem8}, we may think of $I^V\cup I^W$ as a natural
transformation 
\[I^V\cup I^W:h(V\times W)^{\otimes S}(\chi\dim V\times W)
\longrightarrow h(\ol{M})\]
between functors $\gsvw_{\cart}\to(\mbox{graded DMC-motives})$. 

\begin{them} \label{st}
We have
\[I^V\cup I^W=I^{V\times W}.\]
\end{them}
\begin{pf}
This follows from Theorem~\ref{pt} and the identity principle for
DMC-motives, Proposition~8.2 of \cite{BM}.
\end{pf}


\begin{flushleft}
{\tt behrend@math.ubc.ca}
\end{flushleft}

\end{document}